\documentclass[journal]{IEEEtran}
\usepackage{graphicx} %package to manage images
\DeclareGraphicsExtensions{.eps,.pdf,.jpeg,.png}
\usepackage{epstopdf}

\graphicspath{ {Pictures/} }
\ifCLASSINFOpdf
\else
\fi
\usepackage{amsmath,amsfonts,amsthm} % Math packages

\usepackage{enumerate}
\usepackage{siunitx}
\usepackage{eurosym}
\usepackage{booktabs}
\usepackage{accents}
\usepackage{multirow}
\usepackage{ctable}
\usepackage[caption=false,font=footnotesize]{subfig}
\usepackage{fixltx2e}
\usepackage{nicefrac}
\usepackage{algorithm}
\usepackage[noend]{algpseudocode}
\usepackage{color}

\allowdisplaybreaks[1]

\begin{document}
	%
	% paper title
	% can use linebreaks \\ within to get better formatting as desired
	% Do not put math or special symbols in the title.
	\title{{An Online Optimization Algorithm for Alleviating Contingencies in Transmission Networks}}
	%
	%
	% author names and IEEE memberships
	% note positions of commas and nonbreaking spaces ( ~ ) LaTeX will not break
	% a structure at a ~ so this keeps an author's name from being broken across
	% two lines.
	% use \thanks{} to gain access to the first footnote area
	% a separate \thanks must be used for each paragraph as LaTeX2e's \thanks
	% was not built to handle multiple paragraphs
	%

	\author{Nicol\`o~Mazzi,~\IEEEmembership{}%Student~Member,~IEEE,}
		Baosen~Zhang,~\IEEEmembership{Member,~IEEE}
		and~Daniel~S.~Kirschen,~\IEEEmembership{Fellow,~IEEE}% <-this % stops a space
		\thanks{N. Mazzi is with the School of Mathematics, University of Edinburgh, Edinburgh,
			UK (e-mail: nicolo.mazzi@ed.ac.uk).}% <-this % stops a space
		\thanks{B. Zhang and D. S. Kirschen are with the Department
			of Electrical Engineering, University of Washington, Seattle, WA 98195,
			USA (e-mail: zhangbao@uw.edu; kirschen@uw.edu).}% <-this % stops a space
	}
	\maketitle

	% As a general rule, do not put math, special symbols or citations
	% in the abstract or keywords.
	\begin{abstract}
		Power systems are increasingly operated in corrective rather than preventive security mode, which means that appropriate control actions must be taken immediately after a contingency has occurred. This paper proposes an online algorithm for automatically alleviating contingencies such as voltage limit violations and line overloads. Unlike previously proposed approaches, the network itself serves as a natural solver of the power flow equations. This makes it possible to start the implementation immediately and avoids problems caused by modeling errors. Every time the controller receives measurements from the grid, it evaluates the presence of contingencies and computes the optimal corrective actions that can be implemented before the next sampling period, subject to ramping constraints of the generators. These corrective actions are implemented through the standard Automatic Generation Control. Finding the optimal incremental corrective actions is fast because this problem is linearized. The effectiveness of this algorithm at correcting both line overloads and voltage violations is demonstrated using the  IEEE-118 Bus test system.
	\end{abstract}

	% Note that keywords are not normally used for peerreview papers.

	% For peer review papers, you can put extra information on the cover
	% page as needed:
	% \ifCLASSOPTIONpeerreview
	% \begin{center} \bfseries EDICS Category: 3-BBND \end{center}
	% \fi
	%
	% For peerreview papers, this IEEEtran command inserts a page break and
	% creates the second title. It will be ignored for other modes.
	\IEEEpeerreviewmaketitle

	\vspace{-10pt}

	\section{Introduction}
	% The very first letter is a 2 line initial drop letter followed
	% by the rest of the first word in caps.
	%
	% form to use if the first word consists of a single letter:
	% \IEEEPARstart{A}{demo} file is ....
	%
	% form to use if you need the single drop letter followed by
	% normal text (unknown if ever used by IEEE):
	% \IEEEPARstart{A}{}demo file is ....
	%
	% Some journals put the first two words in caps:
	% \IEEEPARstart{T}{his demo} file is ....
	%
	% Here we have the typical use of a "T" for an initial drop letter
	% and "HIS" in caps to complete the first word.
	\IEEEPARstart{T}{raditionally}, power systems have been operated in N-1 preventive security mode, which means that no immediate action is required following a single generation or transmission outage. However, because of the cost of implementing preventive security measures, power systems are increasingly operated in corrective security mode \cite{panciatici2014advanced}. This means that actions must be taken soon after an outage to prevent line overloads from causing cascading outages or voltage violations from leading to a voltage collapse \cite{roald2017corrective}. Various methods have been proposed to calculate what these corrective actions should be \cite{monticelli1987security} and some authors have suggested mechanisms for implementing these measures automatically \cite{bacher1988real}. Typically, these approaches rely on the solution of an Optimal Power Flow (OPF) or Security Constrained Optimal Power Flow (SCOPF) to determine the control actions needed to reach a suitable target operating state. Relying on a model to determine what needs to be done has several drawbacks. First, the model may not accurately represent the behavior of the actual system. Second, the calculation of the target state may require a substantial amount of time and thus delay the implementation of the corrective actions. Third, because the target state cannot be reached instantaneously, other operating constraints might be violated in the process of getting there.

	This paper proposes a closed loop approach to the implementation of post-contingency corrective actions that does not require the solution of the non-linear power flow equations. Instead, incremental control actions based on real-time system measurements are implemented through the Automatic Generation Control (AGC). The network, therefore, acts as a natural solver of the power flow equations, as proposed in \cite{gan2016online}. Corrective actions are thus implemented step by step and take into account what can actually be executed during each AGC cycle. Because they are small and updated in a closed loop, the magnitude and direction of the corrective steps can be calculated using a fast linearized model without causing significant violations of operating constraints.

	The remainder of this paper is organized as follows. Section \ref{sec:literature_review} reviews the relevant literature on corrective actions and closed loop control of power systems. Section \ref{sec:problem_definition} defines the control problem. Section \ref{sec:Contringency Alleviation} describes the optimization models used within the contingency alleviation algorithm. Section \ref{sec:Simulation_Algorithm} presents the algorithms used to simulate the operation of the networks and the contingency alleviation algorithm. Section \ref{sec:test_case} demonstrates the effectiveness of the proposed method using test cases based on the IEEE 118-bus system. Section \ref{sec:Conclusions} concludes.

	\section{Literature Review and Proposed Method}\label{sec:literature_review}

	The approach proposed in this paper relies on three areas of previous work: post-contingency corrective actions, {AGC} and online {OPF} algorithms.

	A significant amount of generation from stochastic renewable energy sources increases the uncertainty that operators must deal with daily. In this context, maintaining power system reliability in the traditional manner {\cite{BaluEtAl1992}} is getting increasingly difficult and costly. Hence the growing interest in post-contingency corrective actions \cite{panciatici2014advanced}. Monticelli et al. \cite{monticelli1987security} are the first to incorporate the possibility of corrective actions in a {SCOPF} and demonstrate that it reduces the operating cost of the system without compromising its operational reliability. Since then, a number of other authors (e.g., \cite{capitanescu2007contingency}, \cite{capitanescu2009coupling}, \cite{chatzivasileiadis2014security}, {\cite{PhanEtAl2014}, and \cite{LiuEtAl2015}}) have developed increasingly sophisticated SCOPF formulations and solution methods to determine the system's optimal operating state when post-contingency corrective actions are possible.

	Other authors have proposed algorithms for determining the optimal set of corrective actions that should be applied when a particular contingency occurs. For example, the authors of \cite{medicherla1981generation} develop an iterative algorithm that modifies the active and reactive power injections at both generator and load buses to alleviate line overloads. Shandilya et al. \cite{shandilya1993method} propose an algorithm for generation rescheduling and load shedding, under the assumption that most of the lines do not operate close to their maximum capacity and that only buses in {the} proximity of the overload are rescheduled. This leads to a non-linear optimization problem that is solved iteratively using a conjugate gradient technique. Bijwe et al. \cite{bijwe1993alleviation} consider both line overloads and voltage violations and use the natural decoupling between active and reactive power to develop two separate models for alleviating line overloads and voltage violations. Similarly, Arini \cite{arini1997fast} proposes an iterative decoupled algorithm for congestion management using linear generation shift distribution factors. Other authors, e.g., \cite{hazarika1998method}, \cite{talukdar2005computationally} and \cite{abbas2016transmission}, develop similar iterative algorithms for congestion management through load shedding and/or generation rescheduling. {References \cite{otomega2007model} and \cite{almassalkhi2015model} propose a model predictive control algorithm to mitigate line overloads}. In recent years fuzzy logic is also used to address this problems {(e.g., \cite{lenoir2009overload} and \cite{pandiarajan2011overload})}.

	The majority of these algorithms aim to evaluate the final values of the control variables (i.e., generally, active power injection and voltage magnitude at generator buses) that would remove the violations of the operating constraints. To do so, they iteratively move the control variables while solving a power flow model to check the status of overloaded lines and bus voltages. However, even when the solution can be found in few seconds, the system may need several minutes to reach the new operating point. Ramp rate constraints on the generators introduce inter-temporal constraints that can not be neglected while considering the implementation of these algorithms. In fact, the solution may not be feasible due to modeling inaccuracies. To avoid these issues, our approach does not solve a power flow to simulate the behavior of the network. Instead, it follows the evolution of its state at each step using actual grid measurement. {When the} measurements are received, it quickly determines the next incremental corrective step by solving a Linear Programming (LP) problem. The inaccuracies introduced by this linearization are not significant because each step is small and deviations are corrected at each step based on actual measurements. This approach also prevents the creation of new violations of operating constraints while correcting the initial one.

	Bacher and van Meeteren \cite{bacher1988real} are the first to propose using real-time corrective actions to track an OPF solution. Gan and Low \cite{gan2016online} recently propose to use the grid as a natural solver of the power flow equations and drastically reduce the computing time required for solving an AC OPF iteratively in radial distribution grids. They claim that, using this approach, it is possible to track the optimal solution when the level of loads or production units changes quickly and continuously. Similar are described in \cite{dall2016design}, \cite{dall2016optimal}, and {\cite{ZhouEtAl2017}}.

\begin{figure}[t!]
	\centering
	\includegraphics[width=0.85\columnwidth]{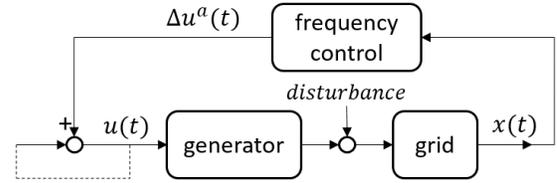}
	\caption{Schematic representation of the frequency control loop. $x$ is the dependent variable (i.e., system's frequency), $u$ is the control variable (i.e., active power production), and $\Delta u^a$ is the variation of $u$ provided by the frequency control in order to compensate the disturbance occurred.}\label{fig.AGC_Control}
\end{figure}

	AGC \cite{JaleeliAGC}, \cite{carpentier1985or} has been used for decades as a closed loop control system to maintain the frequency stability of power systems. The schematic representation shown in Figure \ref{fig.AGC_Control} illustrates how it keeps the frequency (i.e., the output variable $x$) at its nominal level by acting on the active power production of the generators (i.e., the control variables $u$).

\begin{figure}[t!]
	\centering
	\includegraphics[width=0.85\columnwidth]{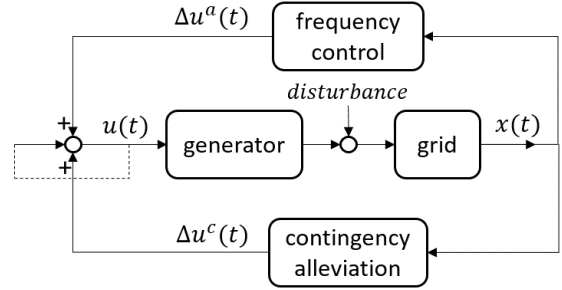}
	\caption{Schematic representation of the integration of the contingency alleviation and frequency controls. $\Delta u^c$ is the adjustments of $u$ (e.g., active power set-points) needed to alleviate the contingency. The output variables $x$ are the system's frequency for the frequency control loop and the flows on the transmission lines for the contingency alleviation loop}\label{fig.Corrective_Control}
\end{figure}

Figure \ref{fig.Corrective_Control} shows how the contingency alleviation concept proposed in this paper could be integrated with the AGC. In this case, the output variables include all the quantities that are measured to detect potential violations of operating limits. The adjustments to the generator active power set-points combine what is needed to maintain the frequency and what is required to alleviate the line overloads. A similar loop adjusts the voltage set-points to correct voltage violations.

	{The idea behind the proposed method is to start adjusting the control variables in the right \emph{direction} almost immediately after a violation of operating constraints has been detected. It will not bring the system to an optimal operating point but will relieve the constraint violations and thus keep it in a acceptable state until the operator moves it towards the solution of an SCOPF. Although based on a linear approximation of the power flow equations, the proposed method is applicable to a wide range of contingencies. For ``soft'' contingencies, our algorithm provides the operators an automatic mechanism to implement corrective actions without having to explicitly initiate these actions. This has the potential to lead to faster correction of minor violations. Our method is also applicable to contingencies that put a heavier stress on the system and make its behavior more non-linear. The linear approximation is still useful as long as the \emph{direction} of the action (e.g., increase or decrease injections) that it suggests is correct. Implementing small steps and continuously using the measurements from the system as feedback takes care of the non-linearity. The proposed approach is thus quite different from solving  the linearized equations to find the best SCOPF as in~\cite{ChiangEtAl2015,ChakrabartiEtAl2014}.}

	{We do not consider the feasibility of real-time adjustments of the protection schemes. However, since the proposed algorithm relieves violations of operating constraints as fast as the ramping limitations of the generators allow, it reduces the likelihood of operation of the protection relays.}

\section{Problem Formulation}\label{sec:problem_definition}

	This section introduces the mathematical model of the network (Section \ref{subsec:Network_Model}) and the formulation of the contingency alleaviation as a control problem (Section \ref{subsec:Control_Problem}).

\subsection{Network Model}\label{subsec:Network_Model}

	We model a meshed power network as a connected graph $G(N^+,E)$ where $N^+:=\{ 0 \} \cup N$, $N:=\{1,2,...,n\} $ and $E \subseteq N^+ \times N^+ $. Each edge of $E$ represents a transmission line and each node of $N^+$ represents a bus. $y_{ij}$, $g_{ij}$ and $b_{ij}$ are respectively the admittance, conductance and susceptance of the transmission line $(i,j) \in E$ ($y_{ij}=g_{ij}+ \textbf{\textit{j}} b_{ij}$). Let be $I_{ij}$ the complex current and $S_{ij} = P_{ij} + \textit{\textbf{j}}Q_{ij}$ the \textit{sending-end} complex power from bus $i$ to bus $j$. The maximum complex power that can be safely transmitted through the line $(i,j)$ is $\overline{S}_{ij}$. A bus $i \in N^+$ can be connected to a generator, a load, both of them or neither. $V_i$ and $\delta_i$ are the voltage magnitude and phase at bus $i \in N^+$. $\overline{v}$ is the maximum acceptable deviation of the voltage magnitude $V_i$ from 1.0 per unit (p.u.). $s_i$, $p_i$ and $q_i$ are the complex, active and reactive power injection at bus $i\in N^+$, respectively.

	Four real variables (or two complex ones) characterize each bus $i \in N$, i.e., $V_i$, $\delta_i$, $p_i$, and $q_i$. Two of these real variables are imposed, while the remaining two are dependent variables determined by the power flow equations. We classify the buses into three categories based on which two variables are imposed, i.e., \textit{slack bus}, \textit{generator bus} and \textit{load bus}. At a \textit{slack bus} $V_i$ and $\delta_i$ are specified, and $p_i$ and $q_i$ are variable. Without loss of generality, we assume that bus 0 is the \textit{slack bus}, assuming for convenience that $V_0 = 1$ p.u. and $\delta_0=0^{\circ}$. For a \textit{generator bus} (PV-bus) $V_i$, $p_i$ are specified, and $\delta_i$ and $q_i$ are variable. For a \textit{load bus} (PQ-bus) $p_i$, $q_i$ are specified, and $V_i$ and $\delta_i$ are variable. Let $N^g$ and $N^l$ be the subsets of PV-buses and PQ-buses, respectively.

	\subsection{Control Problem}\label{subsec:Control_Problem}
	The aim of the paper is to develop an efficient and reliable tool for alleviating contingencies in {transmission} networks. We consider two different types of contingencies, i.e., voltage violations at PQ buses (i.e., $|V_l-1|>\overline{v}$) and line overloads (i.e., $|S_{ij}|>\overline{S}_{ij}$). Starting from an operating state that violates voltage and line flow constraints, we want to determine a series of corrective actions over the control variables $(\mathbf{p}^g,\mathbf{V}^g) = \{(p_g,V_g), g \in N^g \}$ that will bring the system to an operating state that satisfies all operating constraints. Unlike other similar works (e.g., \cite{medicherla1981generation}-\cite{abbas2016transmission}), which iteratively solve an approximate power flow model within the algorithm, we use the measurement of the grid as a solution of the power flow equations. In this way, we avoid the risk of having a solution that is strongly influenced by the approximation used to build the simplified power flow model. Moreover, in order to have a fast response, we only compute the optimal corrective action that can be implemented before the next sampling period, taking ramp rate constraints into account.

	We assume that the measurement from the grid are available every $t^m$ seconds. Each $t^m$, the algorithm receives the values of $\mathbf{V}^l = \{V_l, l \in N^l \}$ and $\mathbf{S} = \{S_{ij}, (ij) \in E \}$, evaluates the presence of violations and modifies the set points $(\mathbf{p}^{set},\mathbf{V}^{set}) = \{(p^{set}_g,V^{set}_g), g \in N^g \}$ in order to reduce the contingency. Moreover, we consider the frequency control response of the generators. In power flow simulations the frequency is assumed to remain constant and the \textit{slack bus} is assumed to produce or absorb whatever power is needed to maintain the load/generation balance. In this work, we model the frequency control as a distributed response of the generators connected to the grid to restore the active power injection $p_0$ at the \textit{slack bus} to its nominal level $\widetilde{p}_0$. We consider that the frequency control intervenes every $t^a$ seconds and modifies the set point of the control variables $\mathbf{p}^g$.
	\begin{figure}[htbp!]
		\centering
		\includegraphics[width=0.95\columnwidth]{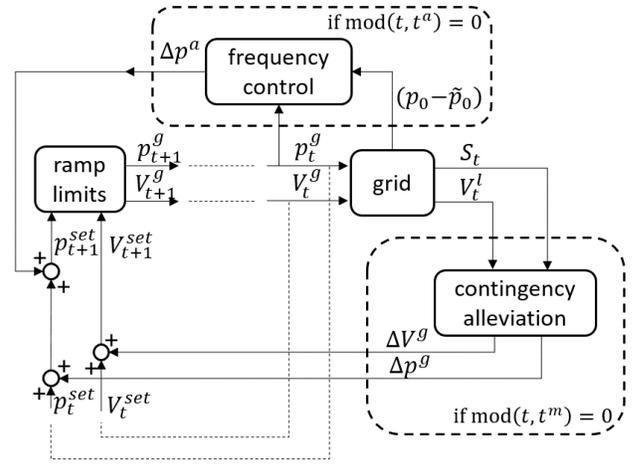}
		\caption{Schematic representation of the control loop of the contingency alleviation algorithm.}\label{fig.Simulation_Blocks}
	\end{figure}

	Figure \ref{fig.Simulation_Blocks} illustrates the control process. For each time step $t$ of the simulation, for given value of $(\mathbf{p}^g,\mathbf{V}^g)$, the \textit{grid} computes the values of the dependent variables $(\mathbf{V}^l,\mathbf{S})$. Every $t^m$ seconds the algorithm determines an optimal step of corrective action to alleviate possible contingencies, as variation of the control variables $(\mathbf{\Delta p}^{g},\mathbf{\Delta V}^{g})$. Similarly, every $t^a$ seconds, the frequency control loop computes the imbalance $p_0-\widetilde{p}_0$ at the slack bus and evaluates $\mathbf{\Delta p}^{a}$ in order to balance it. Then, given the set points $(\mathbf{p}^{set},\mathbf{V}^{set})$, the values of $(\mathbf{p}^g,\mathbf{V}^g)$ are computed for each PV-bus for the next time step $t+1$, considering its ramp rate limits.

	{The proposed method monitors the apparent power flows $S_{ij}$ to detect and quantify line overloads. It could easily be modified to use the line currents $I_{ij}$ for this purpose.}

	\section{Contingency Alleviation}\label{sec:Contringency Alleviation}
	This section defines the contingency alleviation algorithm. Section \ref{subsec:General_Formulation} presents the general formulation of the control problem, while Section \ref{subsec:Linear_Formulation} describes the linearization of the AC power flow equations. Finally, Section \ref{subsec:LP Formulation} formulates the contingency alleviation problem as an LP optimization.

	\subsection{General Formulation}\label{subsec:General_Formulation}
	The first goal of the algorithm is to evaluate the presence of violations of voltage or power flow limits. We initially consider the following objective function:
	\begin{equation}\label{eq:Lmax}
	\mathcal{L} = \mu \sum_{l \in N^{l}} \mbox{max}(|V_l-1| - \overline{v},0) + \dfrac{1}{k}  \sum_{(ij)\in E} \mbox{max}(|S_{ij}|-\overline{S}_{ij},0),
	\end{equation}
	where the parameter $\mu$ controls the weight of the voltage violation term with respect to the line overload one. Parameter $k$ (MVA/pu) is used to convert power flows into per unit. $\mathcal{L}$ is equal to 0 when there are no constraint violations, while it is positive when an operating limit is not respected. {However, the objective function should be sensitive to constraints that are ``close" to be violated, so that further control actions do not cause additional violations of operating constraints. Accordingly, we introduce the function $g(\tau,\epsilon)$ defined as:}
	\begin{equation}\label{eq:g_fun}
	g(\tau,\epsilon) =\begin{cases} 0, & \mbox{if} \hspace{5pt}  \tau \leq - \dfrac{2}{3} \epsilon  \\ \dfrac{1}{3 \epsilon^2} \left( \tau + \dfrac{2}{3} \epsilon \right)^3 , & \mbox{if} \hspace{5pt} - \dfrac{2}{3} \epsilon <  \tau \leq \dfrac{1}{3} \epsilon \\ \tau, &  \mbox{otherwise}.
	\end{cases}
	\end{equation}

\begin{figure}[htbp!]
		\centering
		\includegraphics[width=1.0\columnwidth]{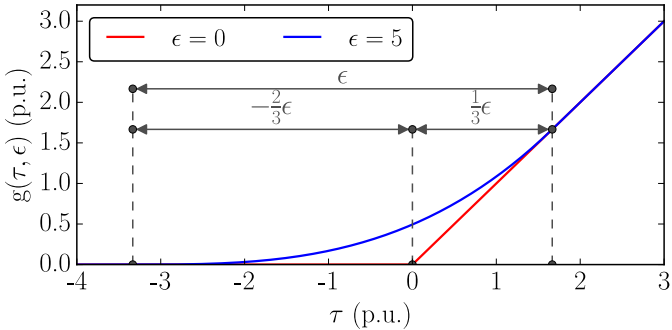}
		\caption{Illustrative example of $g(\tau,\epsilon)$, for $\epsilon=0$ (red) and $\epsilon=5$ (blue).}\label{fig.g_function}
	\end{figure}

	Figure \ref{fig.g_function} illustrates the shape of this function. The red line shows $g(\tau,\epsilon)$ when $\epsilon=0$, which corresponds to the unpenalized function, i.e., $\mathrm{max}(\tau,0)$, while the blue line illustrates $g(\tau,\epsilon)$ when $\epsilon=5$. Note that $g(\tau,\epsilon)\geq 0$ for $\tau \geq - \nicefrac{2}{3} \epsilon$. Therefore, the objective function is greater than 0 when a line or a voltage is close to its limit but there is no actual constraint violation.

	To simplify the notation we introduce the following functions:
	\begin{subequations}
		\begin{align}
		g^V(\tau) =      & g(|\tau-1|-\overline{v}, \xi \; \overline{v} ),\\
		g^S_{ij}(\tau) = & g(|\tau|-\overline{S}_{ij},\xi \; \overline{S}_{ij} ),
		\end{align}
	\end{subequations}
	where $\xi \in [0,1]$ controls the shape of $g^V$ and $g^{S}$. Using these notations, we introduce the penalized modified function:
	\begin{equation}
	\mathcal{L'} = \mu \sum_{l \in N^{l}} g^V(V_l) + \dfrac{1}{k} \sum_{(ij)\in E} g^S_{ij}(S_{ij}).
	\end{equation}

	The aim of the corrective action algorithm is to iteratively modify the control variables $(\mathbf{p}^g,\mathbf{V}^g)$ to alleviate the contingencies. To this effect, we need to determine a quick and reliable path for $(\mathbf{p}^g,\mathbf{V}^g)$ that respects their ramp rate limits. Each $t^m$ seconds, the algorithm receives the measurement of $(\mathbf{V}^l,\mathbf{S})$ and computes the value of $(\mathbf{\Delta p}^g,\mathbf{\Delta V}^g)$ that minimizes $\mathcal{L'}$ while respecting the following ramp rate constraints:
	\begin{subequations}
		\begin{align}
		-R_g t^m \leq \Delta p_g \leq R_g t^m, \; g \in N^g, \label{eq:ramp_lims_p} \\
		-T_g t^m \leq \Delta V_g \leq T_g t^m, \; g \in N^g, \label{eq:ramp_lims_V}
		\end{align}
	\end{subequations}
	where $R_g$ (MW/s) is the ramping limit on $p_g$ and $T_g$ (p.u./s) is the ramping limit on $V_g$. Because voltage regulators are very fast, ramping limits on voltage set-points in \eqref{eq:ramp_lims_V} are not a physical constraints like the ramping constraints on active power set-points in \eqref{eq:ramp_lims_p}. However, we include them because rapid changes in voltage magnitudes are undesirable. If needed, $T_g$ can be set to $\infty$ to remove these constraints.
	At each iteration, we wish to compute the optimal $(\mathbf{\Delta p}^g,\mathbf{\Delta V}^g)$ that can be implemented before the following measurement. We index with the superscript $m$ the values measured from the grid. The optimization problem that we solve at each iteration is:
	\begin{subequations}
		\label{eq:opt0}
		\begin{align}
		\text{Min}
		& \quad \mu \sum_{l \in N^{l}} g^V(V^m_l+\Delta V_l) + \dfrac{1}{k}  \sum_{(ij)\in E} g^S_{ij}(S^m_{ij} + \Delta S_{ij}) \label{eq:opt0_obj}\\
		\text{over}
		& \quad u \; := \; \big(   \Delta p_g, \Delta V_g, \forall g \in N^g \big) \\
		\begin{split}
		& \quad x \; := \; \big( \Delta V_l, \forall l \in N^l; \Delta S_{ij}, \forall (ij) \in E \big)
		\end{split}\\
		\text{s.t.}
		& \quad x = F(u) \label{eq:opt0_PF}\\
		& \quad \sum_{g \in N^g} \Delta p_g = 0\label{eq:opt0_AGC}\\
		& \quad \underline{p}_g \leq p^m_g + \Delta p_g  \leq \overline{p}_g, \quad \forall g \in N^g  \label{eq:opt0_plim}\\
		& \quad |V^m_g + \Delta V_g-1| \leq \overline{v}, \quad \forall g \in N^g  \label{eq:opt0_vlim}\\
		& \quad -R_g t^m \leq \Delta p_g \leq R_g t^m, \quad \forall g \in N^g  \label{eq:opt0_R}\\
		& \quad -T_g t^m \leq \Delta V_g \leq T_g t^m, \quad \forall g \in N^g  \label{eq:opt0_T}
		\end{align}
	\end{subequations}
	where the function $F(\cdot)$ in constraint \eqref{eq:opt0_PF} represents the power flow equations. It computes the variation of the dependent variables $(\mathbf{\Delta V}^l,\mathbf{\Delta S})$, due to $(\mathbf{\Delta p}^g,\mathbf{\Delta V}^g)$. Constraint \eqref{eq:opt0_AGC} imposes the active power balance. Otherwise, eventual imbalances would be compensated by the frequency control loop in the form of a distributed response of the generators. Constraints \eqref{eq:opt0_plim} and \eqref{eq:opt0_vlim} force the active power injection and voltage magnitude at each PV bus to remain within their feasible region. Finally, constraints \eqref{eq:opt0_R} and \eqref{eq:opt0_T} impose limits on the ramp rates. The optimization problem is non-linear due to constraint \eqref{eq:opt0_PF}.

	\subsection{Linear Formulation of $F(\cdot)$}\label{subsec:Linear_Formulation}
	Changes in active and reactive power injections are related to changes in voltage magnitude and phase by:
	\begin{equation}
	\begin{bmatrix} \Delta \mathbf{p} \\ \Delta \mathbf{q}  \end{bmatrix} = - J \begin{bmatrix} \Delta \mathbf{\delta} \\ \Delta \mathbf{V}  \end{bmatrix},
	\end{equation}
	where $J$ is the Jacobian matrix. Factoring $J$ is time consuming for large systems. In order to simplify the computation we use the fast decoupled power flow assumptions \cite{stott1974fast}. We assume that $V_i \approx 1 $ p.u., $\text{cos} \delta_i \approx 1$ and $\text{sin} \delta_i \approx 0$, $  \forall i \in N^+$. Thanks to these assumptions we can decouple the frequency and the voltage control. This leads to:
	\begin{subequations}
		\label{eq:der_FDPF}
		\begin{align}
		&\dfrac{\partial p_i}{\partial \delta_j} \approx - b_{ij} && \quad \dfrac{\partial q_i}{\partial \delta_j} \approx 0,   & \quad i,j \in N^+   \label{eq:fd_delta} \\
		&\dfrac{\partial p_i}{\partial V_j} \approx 0,            && \quad \dfrac{\partial q_i}{\partial V_j} \approx -b_{ij}, & \quad i,j \in N^+    \label{eq:fd_V}
		\end{align}
	\end{subequations}
	The simplified Jacobian is now constant over $t$, so it is sufficient to compute it once at the beginning of the simulation. Aiming in substituting \eqref{eq:opt0_PF} with an alternative linear formulation, we need to evaluate $\frac{\partial \mathbf{S}}{\partial \mathbf{p}^g}$, $\frac{\partial \mathbf{S}}{\partial \mathbf{V}^g}$ and $\frac{\partial \mathbf{V}^l}{\partial \mathbf{V}^g}$.
	For computing $\frac{\partial \mathbf{V}^l}{\partial \mathbf{V}^g}$, we first evaluate the differential $\text{d} \mathbf{q}^l$, where $\mathbf{q}^l = \{q_l, l \in N^l \} $. In accordance with \eqref{eq:fd_V}, we write:
	\begin{equation}
	\text{d}\mathbf{q}^l = \dfrac{\partial \mathbf{q}^l}{\partial \mathbf{V}^l} \text{d}\mathbf{V}^l + \dfrac{\partial \mathbf{q}^l}{\partial \mathbf{V}^g} \text{d}\mathbf{V}^g.
	\end{equation}
	The term $\text{d} \mathbf{q}^l$ is equal to 0 by definition, given that there is no control on the injection of reactive power at a PQ bus. This leads to:
	\begin{equation}
	\dfrac{\partial \mathbf{V}^l}{\partial \mathbf{V}^g} = - \dfrac{\partial \mathbf{q}^l}{\partial \mathbf{V}^g} \left[ \dfrac{\partial \mathbf{q}^l}{\partial \mathbf{V}^l} \right]^{-1}.
	\end{equation}
	Then, computing $\frac{\partial \mathbf{V}}{\partial \mathbf{V}^g}$ is straightforward. Indeed, $\mathbf{V} = \{ V_0 , \mathbf{V}^g , \mathbf{V}^l \}$ and
	\begin{equation}
	\frac{\partial V_0}{\partial V_g} = 0, \;  g \in N^g,
	\end{equation}
	\begin{equation}
	\frac{\partial V_{g'}}{\partial V_g} = \begin{cases} 1, & \mbox{if} \hspace{5pt} g = g' \\
	0, & \mbox{otherwise} \end{cases} \; ,\;  g, g' \in N^g.
	\end{equation}
	To compute the term $\frac{\partial \mathbf{\delta}}{\partial \mathbf{p}^g}$ we invert \eqref{eq:fd_delta}, i.e.,
	\begin{equation}
	\dfrac{\partial \delta_j}{\partial p_g} \approx \left[ \dfrac{\partial \mathbf{p}}{\partial \mathbf{\delta}} \right]^{-1}_{gj}, \quad g \in N^g,\; j \in N^+.
	\end{equation}
	Thanks to the fast decoupled power flow assumptions, we obtain:
	\begin{subequations}
		\label{eq:der_FDPF_PQ}
		\begin{align}
		& \dfrac{\partial P_{ij}}{\partial \mathbf{p}^g} \approx - b_{ij} \left( \dfrac{\partial \delta_i}{\partial \mathbf{p}^g} - \dfrac{\partial \delta_j}{\partial \mathbf{p}^g} \right), && \dfrac{\partial P_{ij}}{\partial \mathbf{V}^g} \approx 0, & (i,j) \in E {}\label{eq:der_FDPF_P} \\
		& \dfrac{\partial Q_{ij}}{\partial \mathbf{V}^g} \approx - b_{ij} \left( \dfrac{\partial V_i}{\partial \mathbf{V}^g} - \dfrac{\partial V_j}{\partial \mathbf{V}^g} \right), && \dfrac{\partial Q_{ij}}{\partial \mathbf{p}^g} \approx 0, & (i,j) \in E
		\end{align}
	\end{subequations}
	Finally, we compute $\dfrac{\partial \mathbf{|S|}}{\partial \mathbf{p}^g}$ and $\dfrac{\partial \mathbf{|S|}}{\partial \mathbf{V}^g}$ as:
	\begin{equation}\label{eq:dS_dpg}
	\dfrac{\partial |S_{ij}|}{\partial \mathbf{p}^g} = \dfrac{P_{ij}}{|S_{ij}|} \; \dfrac{\partial P_{ij}}{\partial \mathbf{p}^g}, \quad (i,j) \in E,
	\end{equation}
	\begin{equation}\label{eq:dS_dVg}
	\dfrac{\partial |S_{ij}|}{\partial \mathbf{V}^g} = \dfrac{Q_{ij}}{|S_{ij}|} \; \dfrac{\partial Q_{ij}}{\partial \mathbf{V}^g}, \quad (i,j) \in E.
	\end{equation}

	\subsection{LP Formulation}\label{subsec:LP Formulation}
	The fast decoupled power flow approximation introduced in Section \ref{subsec:Linear_Formulation} allows us to replace constraint \eqref{eq:opt0_PF} with the following set of constraints:
	\begin{subequations}
		\label{eq:opt01}
		\begin{align}
		& \Delta V_l = \sum_{g \in N^g} \dfrac{\partial V_l}{\partial V_g} \Delta V_g, \quad \forall l \in N^l, \label{eq:opt0_Vl}\\
		& \Delta |S_{ij}| = \sum_{g \in N^g} \left[    \dfrac{\partial |S_{ij}|}{\partial p_g} \Delta p_g +\dfrac{\partial |S_{ij}|}{\partial V_g} \Delta V_g \right], \quad \forall (ij) \in E .  \label{eq:opt0_S}
		\end{align}
	\end{subequations}
	Even though we compute the sensitivities under the decoupled power flow assumptions, the optimization problem \eqref{eq:opt0} can not be decomposed between the control variables $\Delta \mathbf{p}^g$ and $\Delta \mathbf{V}^g$. Since $|S_{ij}| = \sqrt{P^2_{ij} + Q^2_{ij}}$, they both influence $\Delta |S_{ij}|$. However, the linearization of $\partial |S_{ij}|$ in equations \eqref{eq:dS_dpg} and \eqref{eq:dS_dVg} may lead to mistakes when $Q_{ij} \ll P_{ij}$ or, less likely, $P_{ij} \ll Q_{ij}$. Therefore, we prefer to impose some artificial upper limits on $P_{ij}$ and $Q_{ij}$, i.e., $\overline{P}_{ij}$ and $\overline{Q}_{ij}$. We evaluate them by imposing the following conditions:
	\begin{subequations}
		\label{eq:opt02}
		\begin{align}
		& \overline{P}^2_{ij} + \overline{Q}^2_{ij} = \overline{S}_{ij}^2, \quad \forall (ij) \in E  \label{eq:PQS}\\
		& \dfrac{\overline{P}_{ij}}{\overline{Q}_{ij}} = \dfrac{|P_{ij}|}{|Q_{ij}|}, \quad \forall (ij) \in E   \label{eq:PQ_ratio}
		\end{align}
	\end{subequations}
	Thanks to condition \eqref{eq:PQS} we ensure that, if $|P_{ij}| < \overline{P}_{ij}$ and $|Q_{ij}| < \overline{Q}_{ij}$, then $|S_{ij}| < \overline{S}_{ij}$. Then, \eqref{eq:PQ_ratio} imposes that, in case of a violation, the penalty for $|P_{ij}|$ and $|Q_{ij}|$ will be proportional to their contribution to the violation of $\overline{S}_{ij}$. The value of $\overline{P}_{ij}$ and $\overline{Q}_{ij}$ can be computed as
	\begin{equation}\label{eq:Pmax}
	\overline{P}_{ij} = \dfrac{\overline{S}_{ij} }{ \sqrt{ \frac{Q_{ij}^2}{P_{ij}^2} +1 }  }
	\end{equation}
	\begin{equation}\label{eq:Qmax}
	\overline{Q}_{ij} = \dfrac{|Q_{ij}|}{|P_{ij}|} \overline{P}_{ij}
	\end{equation}
	In this way, we can also decompose \eqref{eq:opt0} into two independent subproblems. The first, is solved to evaluate the optimal step $\mathbf{\Delta p}^g$ and its objective function is
	\begin{equation}\label{eq:Lp}
	\mathcal{L}^p = \dfrac{1}{k} \sum_{(ij)\in E} g^P_{ij} (P^m_{ij} + \Delta P_{ij}) + \dfrac{\nu^p}{k} \sum_{g\in N^g} |\Delta p_g|,
	\end{equation}
	where $g^P_{ij}(\cdot)$ is defined as
	\begin{equation}
	g^P_{ij}(\tau) = g(|\tau|-\overline{P}_{ij},\xi \; \overline{P}_{ij} ), \; (ij) \in E.
	\end{equation}
	We introduce the term $\sum_{g\in N^g} |\Delta p_g|$ in \eqref{eq:Lp} to force the model to select only the generators that have a significant capacity to alleviate the contingency. The weight of this term in $\mathcal{L}^p$ is controlled through the parameter $\nu^p$. The optimization problem to evaluate $\mathbf{\Delta p}^g$ is:
	\begin{subequations}
		\label{eq:LP_pg}
		\begin{align}
		{\text{Min}}
		& \quad \mathcal{L}^p   \label{eq:LP_pg_obj}\\
		\text{over}
		& \quad u \; := \; \big(   \Delta p_g,  \forall g \in N^g \big) \\
		& \quad x \; := \; \big( \Delta P_{ij}, \forall (ij) \in E \big) \\
		\text{s.t.}
		& \quad \Delta P_{ij} = \sum_{g \in N^g} \frac{\partial P_{ij}}{\partial p_g} \; \Delta p_g , \; \forall (ij) \in E  \\
		& \quad \sum_{g \in N^g} \Delta p_g = 0 \label{eq:sum0}  \\
		& \quad \underline{p}_g \leq p^m_g + \Delta p_g \leq \overline{p}_g , \; \forall g \in N^g  \\
		& \quad - R_g t^m \leq \Delta p_g \leq R_g t^m , \; \forall g \in N^g
		\end{align}
	\end{subequations}

	The second optimization problem determines the optimal step $\mathbf{\Delta V}^g$, and its objective function is:
	\begin{equation}\label{eq:Lv}
	\begin{split}
	\mathcal{L}^v = & \mu \sum_{l \in N^l} g^{V}(V^m_l(x) + \Delta V_l) +  \\
	&\dfrac{1}{k} \sum_{(ij) \in E} g^{Q}_{ij}(Q^m_{ij}(x)+\Delta Q_{ij}) + \nu^v \sum_{g \in N^g} |\Delta V_g|,
	\end{split}
	\end{equation}
	where $g^Q_{ij}(\cdot)$ is defined as
	\begin{equation}
	g^Q_{ij}(\tau) = g(|\tau|-\overline{Q}_{ij},\xi \; \overline{Q}_{ij} ), \; (ij) \in E.
	\end{equation}
	As in Eq. \eqref{eq:Lp}, we introduce the term $\sum_{g \in N^g} |\Delta V_g|$ and control its weight in \eqref{eq:Lv} through the parameter $\nu^v$. The second optimization problem is:
	\begin{subequations}
		\label{eq:LP_Vg}
		\begin{align}
		\text{Min}
		& \quad  \mathcal{L}^v \label{eq:LP_Vg_obj}\\
		\text{over}
		& \quad u \; := \; \big(   \Delta V_g,  \forall g \in N^g \big) \\
		& \quad x \; := \; \big( \Delta V_l, \forall l \in N^l; \Delta Q_{ij}, \forall (ij) \in E \big) \\
		\text{s.t.}
		& \quad \Delta V_l = \sum_{g \in N^g} \frac{\partial V_l}{\partial V_g} \; \Delta V_g , \; \forall l \in N^l  \\
		& \quad \Delta Q_{ij} = \sum_{g \in N^g} \frac{\partial Q_{ij}}{\partial V_g} \; \Delta V_g , \; \forall (ij) \in E  \\
		& \quad |V^m_g + \Delta V_g -1 | \leq \overline{v} , \; \forall g \in N^g  \\
		& \quad - T_g t^m \leq \Delta V_g \leq T_g t^m , \; \forall g
		\end{align}
	\end{subequations}
	The non-linearities in \eqref{eq:LP_pg} and \eqref{eq:LP_Vg} are replaced with alternative linear formulations. In particular, the functions $g^P(\cdot)$, $g^Q(\cdot)$ and $g^V(\cdot)$ are replaced by piece-wise linear approximations. The absolute values are removed through a conventional linearization technique \cite{boyd2004convex}. The result is an LP model, easily solved using a standard optimization engine \cite{gurobi}.

	\section{Simulation Environment}\label{sec:Simulation_Algorithm}

{This section describes the method used to test the effectiveness of the proposed approach.}
\begin{algorithm}
		\caption{Corrective action}\label{alg:ms}
		\label{alg:corrective_action}
		\textbf{Input:} measurement $(\mathbf{p}^m,\mathbf{V}^m,\mathbf{S}^m)$, line flow limits $\overline{\mathbf{S}}$, maximum voltage deviation $\overline{v}$, relative penalty weights $\nu^p_{\%}$ and $\nu^v_{\%}$.\\
		\textbf{Output:} optimal step $\mathbf{\Delta p}^g$ and $\mathbf{\Delta V}^g$
		\begin{algorithmic}[1]
			\State {compute $\mathcal{L}$ with \eqref{eq:Lmax}}
			\If {$\mathcal{L} > 0$}
			\State compute $\overline{\mathbf{P}}$ and $\overline{\mathbf{Q}}$ with \eqref{eq:Pmax} and \eqref{eq:Qmax}
			\State compute $\mathcal{L}^p$ with \eqref{eq:Lp}, then set $\nu^p = \nu^p_{\%} \mathcal{L}^p $
			\State compute $\mathcal{L}^v$ with \eqref{eq:Lv}, then set $\nu^v = \nu^v_{\%} \mathcal{L}^v $
			\State evaluate $\mathbf{\Delta p^g}$ and $\mathbf{\Delta V^g}$ by solving \eqref{eq:LP_pg} and \eqref{eq:LP_Vg}
			\Else
			\State set $\mathbf{\Delta p^g} = \mathbf{0}$ and $\mathbf{\Delta V^g} = \mathbf{0}$
			\EndIf
		\end{algorithmic}
	\end{algorithm}
	{Algorithm~\ref{alg:corrective_action} is performed by the contingency alleviation block every $t^m$ seconds}. The optimization models \eqref{eq:LP_pg} and \eqref{eq:LP_Vg} are solved only when $\mathcal{L}>0$, where $\mathcal{L}$ defined in Eq. \eqref{eq:Lmax}. Note that $\mathcal{L}$ is greater than 0 if and only if a contingency situation is occurring. {When the reactive power $q_g$ of PV bus $g$ reaches its minimum or maximum limit, it is then modeled as a PQ bus
		where the production of reactive power is fixed at its minimum or maximum limit.}

	{Algorithm~\ref{alg:corrective_action} is integrated in a simulation aimed at reproducing the system's behavior. Static simulations do not consider the system dynamics, but demonstrate the ability of Algorithm~\ref{alg:ms} to alleviate contingencies. Dynamic simulation show how the proposed contingency alleviation algorithm interacts with the system and the frequency control loop.}

	{The proposed algorithm requires as input the network model, measurements of active and reactive power injections ($\mathbf{p}^m$ and $\mathbf{q}^m$) at the PV buses, voltage magnitude $V^m$ at PV and PQ buses, and the power flows $S^m$ in the transmission lines. This is consistent with the data that existing energy management systems can provide using topology processing to evaluate the network model, and the state-estimator for the measurements \cite{almassalkhi2015model}.}
	{\subsection{Static Simulation}}

	\begin{algorithm}
		\caption{Simulation}\label{alg:sim}
		\textbf{Input:} Control variables at $t=0$ $(\mathbf{p}^g_{0},\mathbf{V}^g_{0})$, load consumption $(\mathbf{p}^l,\mathbf{q}^l)$, line flow limits $\overline{\mathbf{S}}$, maximum voltage deviation $v$.
		\begin{algorithmic}[1]
			\For{$t \in \{0,\dots,T\}$}
			\State ($\mathbf{V}^l_t,\mathbf{S}_t) = F(\mathbf{p}^g_t,\mathbf{V}^g_t)$, solved by the network
			\If {$\text{mod}(t,t^m)=0$}
			\State evaluate $\mathbf{\Delta p}^g$ and $\mathbf{\Delta V}^g$ with \textbf{Algorithm~\ref{alg:corrective_action}}
			\State $\textbf{p}^{set} = \mathbf{p}^g_t + \mathbf{\Delta p}^g $
			\State $\textbf{V}^{set} = \mathbf{V}^g_t + \mathbf{\Delta V}^g $
			\EndIf
			\If {$\text{mod}(t,t^a)=0$}
			\State measure imbalance at the slack bus ($p_0-\widetilde{p}_0$)
			\State $\mathbf{\Delta p}^{\text{a}} = - \mathbf{\eta} (p_0-\widetilde{p}_0)  $
			\State $\textbf{p}^{set} = \textbf{p}^{set} + \mathbf{\Delta p}^a $
			\EndIf
			\State $\Delta p_g = \text{max}(\text{min}(p^{set}_g-p_{g,t},R_g),-R_g), \forall g$
			\State $\Delta V_g = \text{max}(\text{min}(V^{set}_g-V_{g,t},T_g),-T_g), \forall g$
			\State $p_{g,t+1} = p_{g,t} + \Delta p_g, \; \forall g$
			\State $V_{g,t+1} = V_{g,t} + \Delta V_g, \; \forall g$
			\EndFor
		\end{algorithmic}
	\end{algorithm}

	{Algorithm~\ref{alg:sim} is used to perform a static simulation of the system's operation.} At each time step $t$, the network solves the power flow equations for given $(\mathbf{p}^g,\mathbf{V}^g)$ and evaluates $(\mathbf{V}^l,\mathbf{S})$. Then, each $t^m$ seconds it performs Algorithm~\ref{alg:corrective_action} to evaluate $\mathbf{\Delta p}^g$ and $\mathbf{\Delta V}^g$ and modifies the set-point of $(\mathbf{ p}^g,\mathbf{ V}^g)$ accordingly. Similarly, each $t^a$ we model the frequency control action by measuring the imbalance at the slack bus and modifying the set-points $(\mathbf{ p}^g,\mathbf{ V}^g)$ in order to compensate for an imbalance. At each time step $t$ we move $(\mathbf{ p}^g,\mathbf{ V}^g)$ towards their set-point $(\mathbf{p}^{set},\mathbf{V}^{set})$, while imposing the ramp rate limits.We simulate the system evolution for $T$ seconds, starting from a contingency.

	{
		In a static environment, a mismatch between power injections and extractions results in a frequency deviation from its nominal value, activating an automatic response from the generators participating in the AGC. In this simplified setup, the frequency is assumed constant, and accordingly the level $p_l$ of the loads do not change during the simulation. However, we approximate frequency control loop as follows. Each $t^a$ seconds, the frequency controller computes the deviation of the active power injection $p_0$ at the slack bus from its nominal value $\widetilde{p}_0$ and provides an immediate distributed response to restore the frequency:
		\begin{equation}
		\Delta p^{\mathrm{a}}_g = - \eta_g \left( p_0 - \widetilde{p}_0 \right), \; \forall g \in N^g,
		\end{equation}
		where $ \eta_g$ define the relative contribution of generator $g$ to the compensation of the system imbalance. The weight $ \eta_g$ is assumed to be:
		\begin{equation}
		\eta_g = \dfrac{\overline{p}_g}{\sum_{g \in N^g} \overline{p}_g}, \; \forall g \in N^g.
		\end{equation}
		It can easily be shown that $\sum_{g} \Delta p^{\textrm{a}}_g + \left( p_0 - \widetilde{p}_0 \right) =  0$.}

{\subsection{Dynamic Simulation}}

	{In a dynamic simulation, the imbalance $(p_0(t)-\widetilde{p}_0)$  between generation and consumption results in the system's frequency $f(t)$ (Hz) deviating from its nominal level $\widetilde{f}$ (i.e., 50 Hz), where $\Delta f(t) = f(t) - \widetilde{f}$. The evolution in time of $\Delta f(t)$ is described by:
		\begin{equation}
		\Delta f(t) = - \dfrac{K_S}{T_S} \int_{0}^{t} \left(p_0(\ell) - \widetilde{p}_0\right) \textrm{d}\ell,
		\end{equation}
		where $K_S$ (Hz/MW) and $T_S$ (s) control the system's dynamic response and $\ell$ is an auxiliary integration variable. The frequency change affects the loads, whose values deviate form their nominal level by $\Delta p^f_l(t)$, equal to:
		\begin{equation}
		\Delta p^f_l(t) = D^f_l \Delta f(t), \quad \forall l \in N^l,
		\end{equation}
		where $D^f_l$ (MW/Hz) is the damping parameter of load $l$. The generators also react to the frequency variation by $\Delta p^g_l(t)$, which can be approximated by:
		\begin{equation}
		\Delta p^f_g(t) = - \dfrac{1}{R^f_g} \Delta f(t), \quad \forall g \in N^G,
		\end{equation}
		where $R^f_g$ (Hz/MW) is the frequency regulation parameter of generator $g$. Finally, the AGC response $\Delta p^{\textrm{a}}_g(t)$ is given by:
		\begin{equation}
		\Delta p^{\textrm{a}}_g (t) = - K^{\textrm{a}}_g \int_{0}^{t} \Delta f(\ell) \hspace{2pt} \textrm{d}\ell, \quad \forall g \in N^g,
		\end{equation}
		where $K^{\textrm{a}}_g$ (MW) controls the AGC response of generator $g$.}

	\section{Test Cases}\label{sec:test_case}

	This section demonstrates the effectiveness of the proposed algorithm on the IEEE 118-Bus system using the parameters listed in Table \ref{tab:params}.
	\begin{table}[htbp!]
		\centering
		\caption{Parameters of the simulation}
		\label{tab:params}
		\begin{tabular}{ccccccccc}
			\specialrule{.1em}{.05em}{.05em}
			$k$ & $\xi $ & $\mu$ & $ \nu^p_{\%} $ & $ \nu^v_{\%} $ & $ R $ & $ T $ & $ t^m $ & $t^a$\\
			(MVA/p.u.)& (-) & (-) & (-) & (-) & (MW/s) & (p.u./s) & (s) & (s) \\
			\midrule
			100 & 0.1 & 5  & 0.001 &  4 & 0.1 & 0.0003 & 4 & 3 \\ \specialrule{.1em}{.05em}{.05em}
		\end{tabular}
	\end{table}

	{\subsection{Static Simulations}\label{subsec:static_simulations}}

	{In the static simulations, three types of contingencies are considered: line overloads, voltage violations, and combinations of line overloads and voltage violations.}\\

\begin{figure}[htbp!]
		\centering
		\includegraphics[width=1.0\columnwidth]{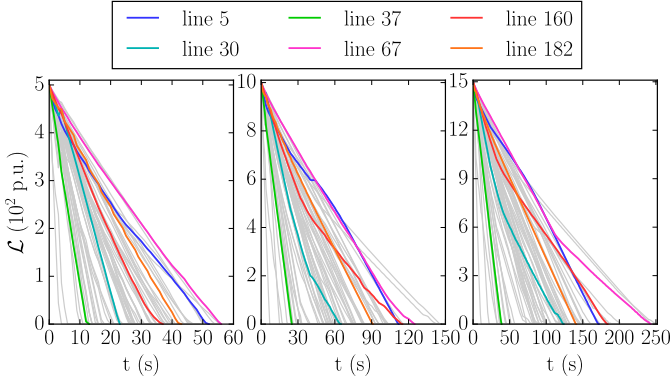}
		\caption{Value of the basic objective function $\mathcal{L}$ as a function of time for overloads $C$ of 5, 10 and 15 MVA on each of the 67 lines considered. Lines 5, 30, 37, 67, 160 and 182 are highlighted to show their behavior during the simulation.}\label{fig.L_3cases}
	\end{figure}

\begin{figure}[htbp!]
		\centering
		\includegraphics[width=0.95\columnwidth]{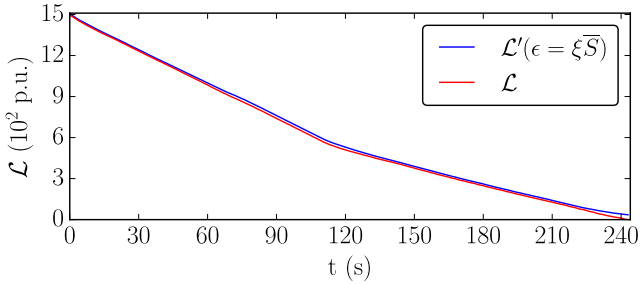}
		\caption{Value of the basic objective function $\mathcal{L}$ (red) and of the penalized objective function $\mathcal{L'}$ (blue) as a function of time for an overload of 15MVA on line 67.}\label{fig.L_L0}
	\end{figure}

\begin{figure}[htbp!]
		\centering
		\includegraphics[width=0.95\columnwidth]{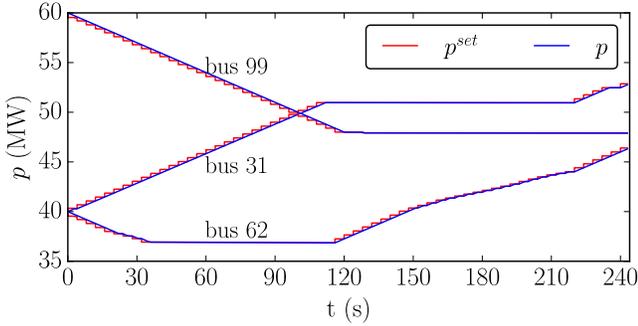}
		\caption{Active power set-points (red)  and actual injections (blue) at buses 30, 61 and 98 to alleviate the overload on line 67.}\label{fig.pg_buses}
	\end{figure}

\begin{figure}[htbp!]
		\centering
		\includegraphics[width=0.95\columnwidth]{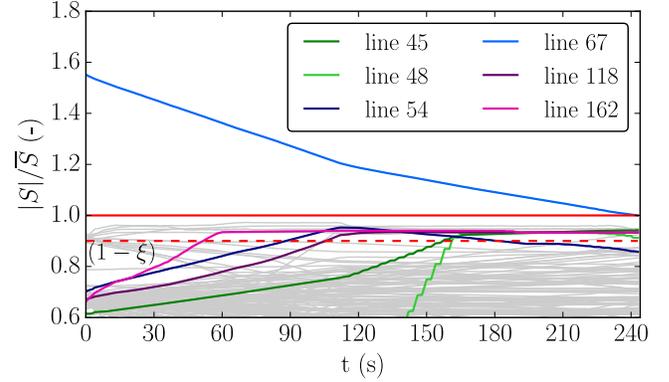}
		\caption{Power flows during the simulation of a 15 MW overload on line 67. Flows on lines 45, 48, 54, 118 and 162 are highlighted to show their behavior during the simulation. All the flows are normalized on the basis of the line rating.}\label{fig.Smtx}
	\end{figure}

	\noindent {\textit{Line Overloads}\\}
	We perform the analysis on 67 lines where the initial flow satisfies all of the following conditions:
	\begin{subequations}
		\label{eq:comp_continuous}
		\begin{align}
		& \dfrac{|P_{ij}|}{|Q_{ij}|} \geq 3  \label{eq:c1}\\
		& |S_{ij}| \geq 20 \text{ MVA}  \label{eq:c2}\\
		& \sum_{g \in N^g} \left| \dfrac{\partial P'_{ij}}{\partial p_g} \right| R_g \geq 0.1 \text{ MW/s} \label{eq:c3}
		\end{align}
	\end{subequations}
	where the derivative $\nicefrac{\partial P'_{ij}}{\partial p_g}$ reflects the effect of the frequency control. Indeed, it considers that if one generator increases its production by 1 MW, all the generators will react to compensate the imbalance of around 1 MW created at the slack bus. Condition \eqref{eq:c3} is imposed to test the algorithm over lines where the system is able to react relatively fast. To create an overload on line $(ij)$, we set:
	\begin{equation}
	\overline{S}_{ij} = S_{ij,t_0} - C
	\end{equation}
	where $C$ represent the magnitude of the overload and $S_{ij,t_0}$ the initial flow on that line. In each simulation we create an overload $C$ of 5, 10 and 15 MVA on only one of these 67 lines. Figure \ref{fig.L_3cases} shows the value of $\mathcal{L}$ over $t$ for each simulation. In all cases, the value of $\mathcal{L}$ is 0 at the end of the simulation, showing that the proposed approach is able to alleviate the contingency in a few minutes.
	Figure \ref{fig.L_L0} compares the evolution of $\mathcal{L}$ and $\mathcal{L'}$ over time for an overload of 15 MW on line 67. $\mathcal{L}$ decreases at a constant rate up to around $t=110$ s, after which point it decreases at a slower rate. The explanation for this change can be found on Figs. \ref{fig.pg_buses} and \ref{fig.Smtx}. Figure \ref{fig.pg_buses} shows the changes in active power injections that the algorithm implements to remove this overload. Figure \ref{fig.Smtx} shows how the flows in all the lines change as a result of these corrective actions. The area between the continuous and dashed red lines is the region where the penalized objective function is sensitive to flows that are close to their maximum. At $t \approx$  110 s, line 118 enters in this penalized region and forces a change in the way the algorithm adjusts the control variables. At $t \approx$  240 s, $\mathcal{L} = 0$ indicating that the overload has been removed. At that time $\mathcal{L'}$ is not quite zero because some line flows are close to their limit.\\

	{We also tested how the corrective action algorithm performs when inaccurate line susceptances are used to compute the sensitivities, e.g., $\frac{\partial \mathbf{S}}{\partial \mathbf{p}^g}$. Let $b_{ij}$ be the actual susceptance of line $(ij)$ and $\widetilde{b}_{ij}$ the value used in the algorithm. $\widetilde{b}_{ij}$ is given by:
		\begin{equation}
		\widetilde{b}_{ij} = \beta_{ij} b_{ij}, \quad \forall (ij) \in E,
		\end{equation}
		where $\beta_{ij} \sim \mathcal{N} \left( \mu_{\beta}, \sigma^2_{\beta} \right)$ and $\mu_{\beta}=1$, $\sigma_{\beta}=0.2$. This translates in an expected difference between $\widetilde{b}_{ij}$ and $b_{ij}$ of around 16\%. Parameter $\beta_{ij}$ is different for each susceptance and is added to all the susceptances simultaneously.
		We ran 200 simulations for a contingency $C$ of 10 MW on lines 5, 30, 37, 67, 160, and 182. These simulations also consider that some measurements (e.g., $S_{ij}$) may not always be updated in time. To simulate this effect, at each iteration, the power flow in one randomly selected line is not updated. The algorithm uses the previous value for this flow measurement. Table \ref{tab:sims2} summarizes the results and shows that the proposed algorithm is still able to remove the violations in 100\% of the tests performed.}

	\begin{table}[htbp!]
		\centering
		{
			\caption{Simulations with errors in the estimated lines susceptance and communication failures.}
			\label{tab:sims2}
			\begin{tabular}{cccccc}
				\specialrule{.1em}{.05em}{.05em}
				{line} & $C$ & \# of tests  & \multicolumn{3}{c}{elapsed time} \\
				&      &      & average & min & max \\
				& (MW) & (-) & (s) & (s) & (s) \\
				\midrule
				5   & 10 & 200  & 115.4 &  101 & 137 \\
				30  & 10 & 200  & 63.3 &  57 & 67 \\
				37  & 10 & 200  & 25.5 &  25 & 27 \\
				67  & 10 & 200  & 125.2 &  122 & 135 \\
				160 & 10 & 200  & 113.7 &  107 & 133 \\
				182 & 10 & 200  & 91.0 &  91 & 91 \\ \specialrule{.1em}{.05em}{.05em}
			\end{tabular}}
		\end{table}

		\noindent {\textit{Voltage Violations}\\}

\begin{figure}[htbp!]
			\centering
			\includegraphics[width=0.95\columnwidth]{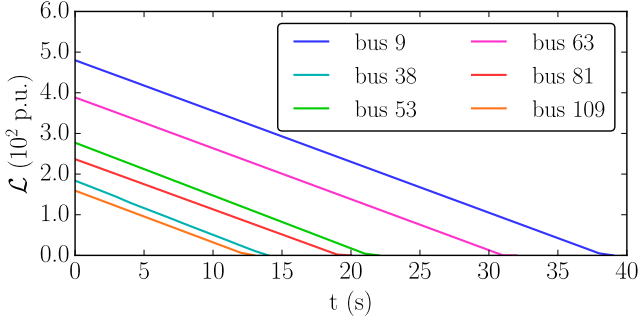}
			\caption{Evolution of the basic objective function $\mathcal{L}$ for voltage violations at buses 9, 38, 53, 63, 81 and 109.}\label{fig.Lvv_buses}
		\end{figure}

\begin{figure}[htbp!]
			\centering
			\includegraphics[width=0.95\columnwidth]{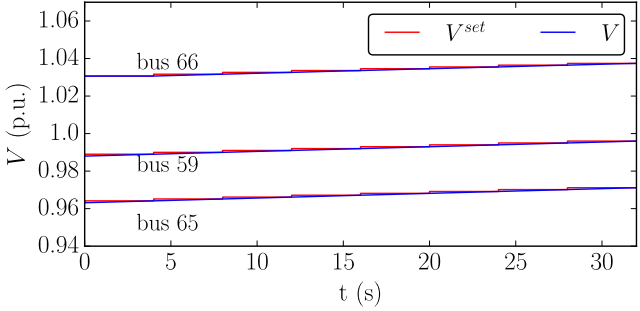}
			\caption{Voltage set-points (red) and actual voltages (blue) at PV buses 59, 65 and 66 for the under-voltage contingency at bus 63.}\label{fig.Vg_buses}
		\end{figure}

\begin{figure}[htbp!]
			\centering
			\includegraphics[width=0.95\columnwidth]{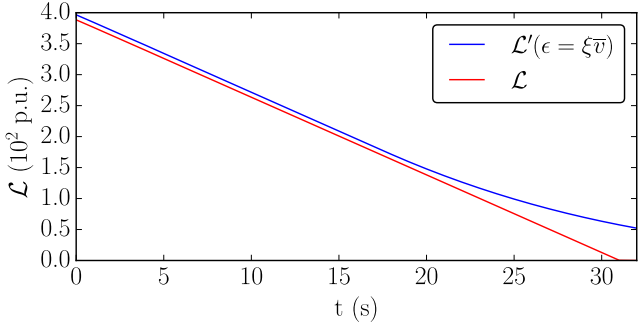}
			\caption{Evolution of the basic $\mathcal{L}$ (red) and penalized $\mathcal{L'}$ (blue) objective functions for a voltage violation at bus 63.}\label{fig.L_L0vv}
		\end{figure}

\begin{figure}[htbp!]
			\centering
			\includegraphics[width=0.95\columnwidth]{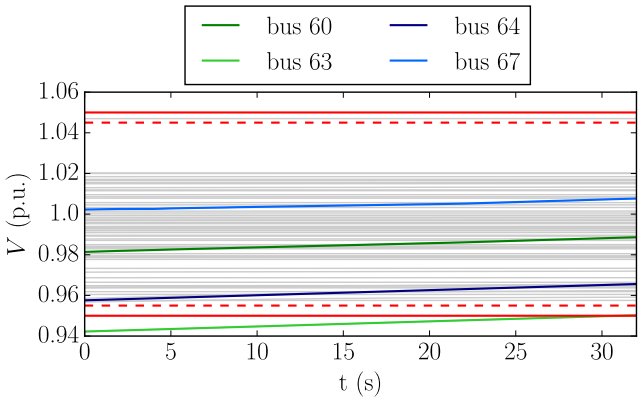}
			\caption{Voltage at the PQ buses for the under-voltage contingency at bus 63. Buses 60, 64 and 67 are highlighted to show their behavior during the simulation. }\label{fig.Vlmtx}
		\end{figure}

		To create voltage violations we modified the reactive power injection separately at PQ buses 9, 38, 53, 63, 81 and 109.
		%Figure \ref{fig.V_contingency} shows the voltage magnitude at all buses except the slack bus excluded before these changes in injection ($t=0$). The blue and gray circles show the voltages at the PV and PQ buses before these changes in reactive power injections. The red empty circles identify the buses where we create the voltage violations, while the red filled circles show the voltages after these changes in reactive injection.
		%
		Figure \ref{fig.Lvv_buses} shows how the value of $\mathcal{L}$ evolves in response to the corrective actions deployed for these 6 voltage violations. Figure \ref{fig.Vg_buses} shows how the algorithm adjusts the voltage set-points at PV buses 59, 65 and 66 to correct a voltage violation at bus 63. Figure \ref{fig.L_L0vv} compares the values of $\mathcal{L}$ and $\mathcal{L'}$ over $t$ for this voltage violation at bus 63. Figure \ref{fig.Vlmtx} illustrates the evolution of the voltage magnitudes at all PQ buses. Most of the voltages are not affected by the corrective actions because voltage problems are typically local. At the end of the simulation all the voltages are within the acceptable range, including bus 63.\\

		% \begin{figure}[htbp!]
		%   \centering
		%   \includegraphics[width=0.7\columnwidth]{V_contingency}
		%  \caption{Voltage magnitudes at the PQ (gray) and PV (blue) buses before and after (red) the changes in reactive injections.}\label{fig.V_contingency}
		% \end{figure}

		\noindent {\textit{Simultaneous Voltage Violations and Line Overloads}\\}

\begin{figure}[htbp!]
			\centering
			\includegraphics[width=0.95\columnwidth]{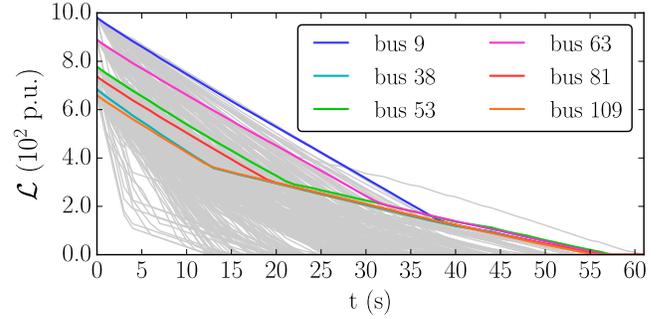}
			\caption{{$\mathcal{L}$ over $t$ for an overload of 5 MVA on line 67 and voltage violations at bus 9, 38, 53, 63, 81 and 109.}}\label{fig.L_buses}
		\end{figure}

		{The last set of static simulations considers simultaneous voltage violation and line overload. We ran the simulations considering each combination of the voltage violations and the line overloads of Section \ref{subsec:static_simulations}}, for C of 5 MVA. Figure \ref{fig.L_buses} illustrates how $\mathcal{L}$ evolves during these 402 simulations, all of which converge in less than one minute. In this figure, the 6 simulations involving an overload of line 67 are highlighted. In all these 6 cases, the rate at which $\mathcal{L}$ decreases at one point during the simulation. This is because voltage violations are corrected more quickly than line overloads because the active power output of generators changes more slowly than their terminal voltage. This change in rate occurs when the voltage violation has been cleared, while line 67 is still overloaded.

		{\subsection{Dynamic Simulations}\label{subsec:dynamic_simulations}}

		{In the dynamic simulation, the parameters $K_S$ and $T_S$ are set at 0.05 Hz/MW and 10 s, respectively. The damping parameter of the loads $D^f_l$ is:
			\begin{equation}
			D^f_l = D^f \dfrac{p_l(0)}{\sum_l p_l(0)}, \quad \forall l \in N^l,
			\end{equation}
			where $p_l(0)$ is the nominal value of the load  and $D^f$ is 20 MW/Hz. $R^f_g$ and $K^{\textrm{a}}_g$ are given by:
			\begin{equation}
			R^f_g = R^f \dfrac{\sum_g \overline{p}_g}{ \overline{p}_g}, \quad \forall g \in N^g,
			\end{equation}
			\begin{equation}
			K^{\textrm{a}}_g = K^{\textrm{a}} \dfrac{\overline{p}_g}{ \sum_g \overline{p}_g}, \quad \forall g \in N^g,
			\end{equation}
			where $R^f$ and $K^{\textrm{a}}$ are 0.005 Hz/Mw and 80 MW, respectively.}

\begin{figure}[htbp!]
			\centering
			\includegraphics[width=0.95\columnwidth]{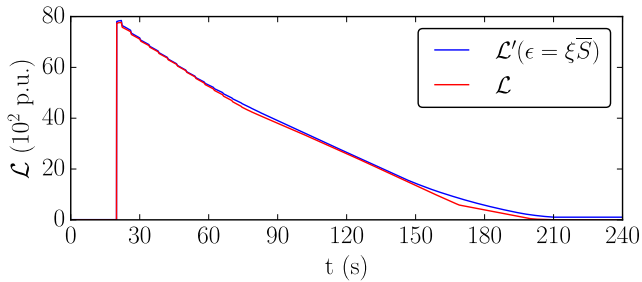}
			\caption{{Evolution of the basic $\mathcal{L}$ (red) and penalized $\mathcal{L'}$ (blue) objective functions for a failure of line 30 at $t=20$ s.}}\label{fig.L_L0_mult}
		\end{figure}

\begin{figure}[htbp!]
			\centering
			\includegraphics[width=0.95\columnwidth]{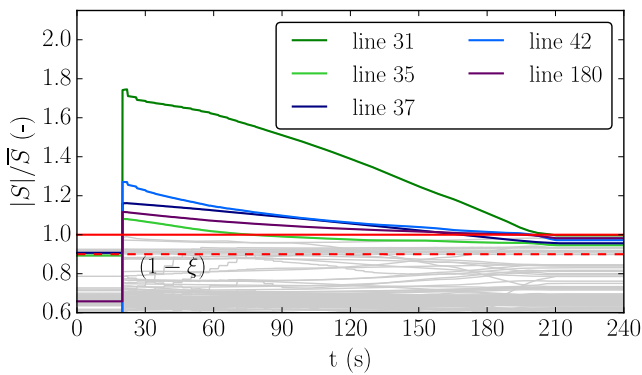}
			\caption{{Power flows during the simulation of a  failure of line 30 at $t=20$ s. Flows on lines 31, 35, 37, 42 and 180 are overloaded after the contingency. All the flows are normalized on the basis of the line rating.}}\label{fig.Smtx_mult}
		\end{figure}

\begin{figure}[htbp!]
			\centering
			\includegraphics[width=0.95\columnwidth]{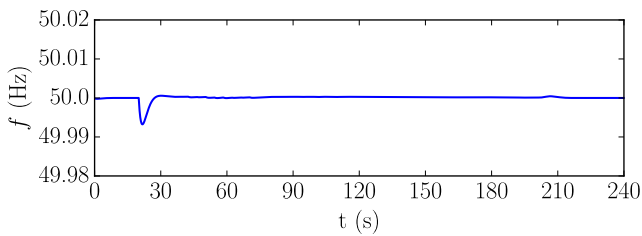}
			\caption{{System's frequency $f(t)$ during the simulation of a  failure of line 30 at $t=20$ s.}}\label{fig.freq_mult}
		\end{figure}

		{
			We use the dynamic simulation to analyze the consequences of the outage of line 30 at $t=20$ s, which causes overloads on lines 31, 35, 37, 42, and 180. In this case the algorithm needs to recompute the sensitivity matrices because the topology of the grid has changed. Calculating these sensitivities requires a one-time computation of about 80 ms, while the LP problems \eqref{eq:LP_pg} and \eqref{eq:LP_Vg} are solved in around 20 ms (each $t^m$ seconds).
			Figure \ref{fig.L_L0_mult} shows how the value of $\mathcal{L}$ evolves as the corrective actions are implemented by the proposed method. $\mathcal{L}$ is 0 up to 20 s, as no violation has yet arisen. After the failure of line 30, $\mathcal{L}$ increases to almost 8000 p.u. and  is then brought back to 0 (i.e., no overloads) in around 180 s. Figure \ref{fig.Smtx_mult} illustrates the normalized flows in the transmission lines during this dynamic simulation, highlighting the 5 lines that overloaded after the failure of line 30. In particular, line 31 is severely overloaded (almost 80\% above its maximum limit) and is  the last one to be returned to an acceptable value. Figure \ref{fig.freq_mult} shows how the system's frequency $f$ evolves during the simulation. Note  the small drop in frequency at $t=20$ s (i.e., when line 30 fails), which is rapidly restored by the AGC. After this event, the frequency is stable and barely affected by the corrective actions because the proposed algorithm coordinates the adjustments in the power outputs of the generator buses in a way that ensures that the total power injection is constant (see constraint \eqref{eq:sum0} in model \eqref{eq:LP_pg}).}

\section{Conclusion}\label{sec:Conclusions}

		This paper presents an innovative algorithm to alleviate contingencies in {transmission} networks. Unlike other techniques described in the literature, the proposed approach uses the network as a natural solver of the power flow equations. This method can be used in parallel with the frequency control loop because they would work on a similar time scale. Every time a measurement from the grid is available, the algorithm evaluates the optimal step in the control variables required to alleviate a line overload, a voltage violation or a combination of the two. Given that the ramp rate constraints on the generators prevent large changes between two consecutive sampling periods, assuming a linear behavior of the system with respect to the control variables is a very good approximation. This allows us to formulate an LP optimization problem, which can be implemented efficiently even for large networks. The effectiveness of this approach is tested on the IEEE 118-Bus system, for line overloads, voltage violations and combinations of these two types of contingencies. These tests demonstrate that it is able to alleviate contingencies quickly, without causing violations of other operating constraints.

		\ifCLASSOPTIONcaptionsoff
		\newpage
		\fi

		% trigger a \newpage just before the given reference
		% number - used to balance the columns on the last page
		% adjust value as needed - may need to be readjusted if
		% the document is modified later
		%\IEEEtriggeratref{8}
		% The "triggered" command can be changed if desired:
		%\IEEEtriggercmd{\enlargethispage{-5in}}

		% references section

		\bibliographystyle{IEEEtran}
		\bibliography{bibliography}

\vfill

		% insert where needed to balance the two columns on the last page with
		% biographies
		%\newpage

		%\begin{IEEEbiographynophoto}{Jane Doe}
		%Biography text here.
		%\end{IEEEbiographynophoto}

		% You can push biographies down or up by placing
		% a \vfill before or after them. The appropriate
		% use of \vfill depends on what kind of text is
		% on the last page and whether or not the columns
		% are being equalized.

		%\vfill

		% Can be used to pull up biographies so that the bottom of the last one
		% is flush with the other column.
		%\enlargethispage{-5in}

		% that's all folks
\end{document}